\documentclass{article}

\usepackage{graphicx}

\usepackage{amsmath}
\usepackage{amssymb}
\usepackage{indentfirst}
\usepackage{enumerate}

\newtheorem{thm}{Theorem}[section]
\newtheorem{cor}[thm]{Corollary}
\newtheorem{lem}[thm]{Lemma}
\newtheorem{de}[thm]{Definition}

\newcounter{bean}
\newcounter{milk}

\begin{document}

\title{The Mystery of the Shape Parameter}         % Enter your title between curly braces
\author{Lin-Tian Luh\\Department of Financial and Computational Mathematics, Providence University\\Shalu, Taichung, Taiwan\\
Email:ltluh@pu.edu.tw\\Fax:886-4-26324653, Tel:886-4-26328001ext.15126 }        % Enter your name between curly braces
\date{\today}          % Enter your date or \today between curly braces
\maketitle
{\bf Abstract}. In this paper we deal with the choice of the  shape parameter $c$ contained in the multiquadrics $(-1)^{\lceil \frac{\beta}{2} \rceil}(c^{2}+\|x\|^{2})^{\frac{\beta}{2}}$ . The optimal value of $c$ was unknown and has obsessed experts in the field of radial basis functions (RBFs) since 1971. People conversant with RBFs have developed various algorithms or models to overcome this question due to its extreme importance when applying RBFs to a variety of scientific fields, including numerical interpolations, PDEs, molecular quantum mechanics, etc. However, hitherto, no satisfactory criteria are available, both in theory and practice. Lack of theoretical understanding to this question probably is the main trouble. The purpose of this paper is to uncover its mystery and present a practically useful method to choose it.

Although the author has in recent years published a few papers dealing with this question, they all are based on evenly spaced data points. Here we deal with its choice in a purely scattered data setting.

\section{Introduction}
We first clarify the definition of the radial function we are concerned about. In the Abstract, $\lceil \frac{\beta}{2} \rceil$ means the smallest integer greater than or equal to $\frac{\beta}{2}$. In order to cite a core theorem of Luh \cite{Lu1}, a modified definition of the radial function is adopted here. In this paper the radial function $h$ is defined by the formula

\begin{align}
h(x):=\Gamma (-\frac{\beta }{2})(c^{2}+\| x\|^{2})^{\frac{\beta }{2}},\quad \beta \in R\setminus 2N_{\geq 0},\quad c>
0,\quad x \in R^{n},
\end{align}
where $\Gamma(\cdot)$ is the classical gamma function which is used only to control the sign of $h(x)$. This function is very powerful because it can be used to make interpolation for scattered data points. The interpolant is of the form

\begin{equation}
s(x):=\sum_{i=1}^{N}c_{i}h(x-x_{i})+p(x)
\end{equation}
where $p(x)\in P_{m-1}$, the space of polynomials of degree less than or equal to $m-1$ in $R^{n}$, $X=\{ x_{1}, \cdots ,x_{N}\}$ is the set of interpolation centers. For $m=0,\ P_{m-1}:=\{0\}$. For any given function $f(\cdot)$, we require that $s(\cdot )$ interpolate $f(\cdot )$ at data points $(x_{1},f(x_{1})),\cdots , (x_{N},f(x_{N}))$. Therefore a linear system constituted by
\begin{equation}
\sum_{i=1}^{N}c_{i}h(x_{j}-x_{i})+\sum_{i=1}^{Q}b_{i}p_{i}(x_{j})=f(x_{j}),\ j=1,\cdots,  N
\end{equation}
and 
\begin{equation}
\sum_{i=1}^{N}c_{i}p_{j}(x_{i})=0 \hspace{3.5cm} ,\ j=1, \cdots,  Q,
\end{equation}
where $\{ p_{1},\cdots , p_{Q}\}$ is a basis of $P_{m-1}$, has to be satisfied. Since $h$ is c.p.d., this requirement will be theoretically satisfied, as shown in Madych et al. \cite{MN1} and Wendland \cite{We}. However if $c$ is very large, $h$ will be numerically constant, making the linear system (3),(4) numerically unsolvable. Moreover, if $c$ is very large, the coefficient matrix of the linear system will have a very large condition number, making the interpolating function $s(\cdot)$ unreliable when $f(x_{1}),\cdots ,f(x_{N})$ are not accurately evaluated, as pointed out by Madych in \cite{MN2}.

\subsection{Fundamental theory}
Our approach of choosing $c$ is based on Theorem 2.4 and Corollary 2.5 of Luh \cite{Lu1} which we will take directly but make a slight modification to make them easier to understand.

Before introducing the main theorem, we need some fundamental definitions. Let ${\cal D}(R^{n})$ denote the space of complex-valued functions on $R^{n}$ that are compactly supported and infinitely differentiable. For any integer $m\geq 1$, let
$${\cal D}_{m}=\{\varphi \in {\cal D}(R^{n}):\ \int x^{\alpha}\varphi (x)dx=0\ for\ all\ |\alpha|<m\}.$$
If $m=0,\ {\cal D}_{m}:={\cal D}(R^{n})$.
\begin{de}
Let $h$ be as in (1) and $m=max\{ 0, \lceil \frac{\beta}{2} \rceil \}$. We write $f\in {\cal C}_{h,m}(R^{n})$ if $f\in C(R^{n})$ and there is a constant $c(f)$ such that for all $\varphi$ in ${\cal D}_{m}$,
\begin{equation}
\left| \int f(x)\varphi(x)dx\right| \leq c(f)\left\{ \int \int h(x-y)\varphi(x)\overline{\varphi(y)}dxdy\right\} ^{1/2}.
\end{equation}
If $f\in {\cal C}_{h,m}(R^{n})$, we let $\| f\| _{h}$ to denote the smallest constant $c(f)$ for which (5) is true.
\end{de}

The function space ${\cal C}_{h,m}(R^{n})$, abbreviated as ${\cal C}_{h,m}$, is called the native space whose characterization can be found in \cite{Lu2, Lu3,  MN1, MN3, We}.   
\begin{de}
For $n=1,2,3,\cdots $, the sequence of integers $\gamma_{n}$ is defined by $\gamma_{1}=2$ and $\gamma_{n}=2n(1+\gamma_{n-1})$ if $n>1$.
\end{de}
\begin{de}
Let $n$ and $\beta$ be as in (1). The numbers $\rho$ and $\Delta_{0}$ are defined as follows.
\begin{list}
  {(\alph{bean})}{\usecounter{bean} \setlength{\rightmargin}{\leftmargin}}
  \item Suppose $\beta <n-3$. Let $s=\lceil \frac{n-\beta -3}{2}\rceil $. Then 
    \begin{list}{(\roman{milk})}{\usecounter{milk} \setlength{\rightmargin}{\leftmargin}}
      \item if $\beta <0,\ \rho=\frac{3+s}{3}\ and\  \Delta_{0}=\frac{(2+s)(1+s)\cdots 3}{
 \rho^{2}};$
      \item if $\beta >0,\ \rho=1+\frac{s}{2\lceil \frac{\beta}{2}\rceil +3} \ and \ \Delta_{0}=\frac{(2m+2+s)(2m+1+s)\cdots (2m+3)}{\rho^{2m+2}}$ \\
where $ m=\lceil \frac{\beta}{2}\rceil$.          
    \end{list}
  \item Suppose $n-3\leq \beta <n-1$. Then $\rho=1$ and $\Delta_{0}=1$.
  \item Suppose $\beta \geq n-1$. Let $s=-\lceil \frac{n-\beta -3}{2}\rceil $. Then
 $$\rho =1\ and \ \Delta_{0}=\frac{1}{(2m+2)(2m+1)\cdots (2m-s+3)} \ where \ m=\left\lceil \frac{\beta}{2}\right\rceil.$$  
\end{list}
\end{de}

We can now introduce our main theorem. 

\begin{thm}
  Let $h$ be defined as in (1) and $m=max\{ 0, \lceil \frac{\beta}{2}\rceil \}$. Then given any positive number $b_{0}$, there are positive constants $\delta_{0}$ and $\lambda,\ 0<\lambda <1$, which depend completely on $b_{0}$ and $h$,  for which the following is true: For any cube $E$ in $R^{n}$ of side length $b_{0}$, if $f\in {\cal C}_{h,m}$ and $s$ is the map defined as in (2) which interpolates $f$ on a finite subset $X$ of $E$, then
\begin{equation}
  |f(x)-s(x)| \leq 2^{\frac{n+\beta+1}{4}}\pi^{\frac{n+1}{4}}\sqrt{n\alpha_{n}}c^{\frac{\beta}{2}}\sqrt{\Delta_{0}}(\lambda)^{\frac{1}{\delta}}\| f\| _{h}
\end{equation} 
holds for all $0<\delta \leq \delta_{0}$ and all $x$ in $E$ provided that $\delta=d(E,X):=\sup _{y\in E}\inf _{x\in X}|y-x|$. Here, $\alpha_{n}$ denotes the volume of the unit ball in $R^{n}$, and $c,\ \Delta_{0}$ were defined in (1) and Definition 1.3 respectively. Moreover $\delta_{0}=\frac{1}{6C\gamma_{n}(m+1)}$, and $\lambda=(\frac{2}{3})^{\frac{1}{6C\gamma_{n}}}$ where
$$C=\max \left\{ 2\rho'\sqrt{n}e^{2n\gamma_{n}},\ \frac{2}{3b_{0}}\right\},\ \rho'=\frac{\rho}{c}.$$
The integer $\gamma_{n}$ was defined in Definition 1.2, and $\| f\| _{h}$ is the $h$-norm of $f$ in ${\cal C}_{h,m}$, as defined in Definition 1.1. The constant $\rho$ was defined in Definition 1.3.

\end{thm}
{\bf Remark}: Theorem 1.4 was taken directly from Luh \cite{Lu1} with only a slight modification. The proof is very technical and absolutely nontrivial. The main contribution of the theorem is that it unveils the mystery of $\lambda$ and $\delta_{0}$ whose values were unknown. Both numbers appear in the exponential-type error bound for multiquadric interpolation presented in Madych et al. \cite{MN3} which is only an existence theorem. The clarification of these constants had been considered to be a hard question. Obviously the domain $E$ in Theorem 1.4 can be extended to a more general set $\Omega \subseteq R^{n}$ which can be expressed as the union of rotations and translations of a fixed cube of side $b_{0}$.
In fact, Madych presented two kinds of exponential-type error bounds altogether, one in \cite{MN2}, the other in \cite{MN3}. The former was used by Madych in \cite{MN2} to establish a criterion of choosing $c$ which is very restrictive. According to his criterion, $c$ is required to be in the interval $[0,\ b_{0}]$ where $b_{0}$ denotes the side length of the domain cube. However, the optimal choice of $c$ may be much larger than $b_{0}$ and cannot be tested by his approach. In this paper we use the error bound in \cite{MN3} to handle it. This approach proves to be very successful. The value of $c$ can be arbitrarily large.

In (6) it is clearly seen that the error bound is greatly influenced by the shape parameter $c$. However, in order to present useful criteria for the choice of $c$, we still need some work to make things transparent.
\begin{de}
  For any $\sigma >0$, the class of band-limited functions $f$ in $L^{2}(R^{n})$ is defined by
$$B_{\sigma}=\{ f\in L^{2}(R^{n}):\ \hat{f}(\xi)=0\ if\ |\xi|>\sigma\},$$
where $\hat{f}$ denotes the Fourier transform of $f$.
\end{de}
\begin{lem}
  Let $h$ be as in (1) with $\beta>0$. Any function $f$ in $B_{\sigma}$ belongs to ${\cal C}_{h,m}$ and 
\begin{equation}
  \| f\| _{h}\leq \sqrt{m!S(m,n)}2^{-n-\frac{1+\beta}{4}}\pi^{-n-\frac{1}{4}}\sigma^{\frac{1+\beta+n}{4}}e^{\frac{c\sigma}{2}}c^{\frac{1-\beta-n}{4}}\| f\| _{L^{2}(R^{n})},
\end{equation}
where $c,\ \beta$ are as in (1) and $S(m,n)$ is a constant determined by $m$ and $n$.
\end{lem}
{\bf Proof}. We are going to show $B_{\sigma}\subseteq {\cal C}_{h,m}$ by Corollary 3.3 and Theorem 5.2 of Madych et al. \cite{MN3}.

Let $m=\lceil \frac{\beta}{2}\rceil$. Then the requirement $a_{\gamma}=0$ for all $|\gamma|=2m$ in Corollary 3.3 of \cite{MN3} is an immediate result of Theorem 5.2 of \cite{MN3}. Since $B_{\sigma}\subseteq L^{2}(R^{n})$, any $f\in B_{\sigma}$ is a member of ${\cal S}'$ where ${\cal S}$ denotes the Schwarz space. Now, let $\rho(\xi)$ be the Borel measure mentioned in Corollary 3.3 of \cite{MN3}. By (3.9) of \cite{MN3}, it suffices to show that
$$\| f\| _{h}:=\left\{ \sum_{|\alpha|=m}\left(\frac{m!}{\alpha!}\right)\| (D^{\alpha}f)^{\hat{}}\| ^{2}_{L^{2}(\rho)}\right\} ^{1/2}<\infty$$
for all $f\in B_{\sigma}$. We proceed as follows.
\begin{eqnarray*}
  &   & \left\{ \sum_{|\alpha|=m}\frac{m!}{\alpha!}\int _{R^{n}}| (D^{\alpha}f)^{\hat{}}(\xi)|^{2}d\rho(\xi)\right\} ^{1/2} \\
  & = & \left\{\sum _{|\alpha|=m}\frac{m!}{\alpha!}\int _{R^{n}}|i^{m}\xi^{\alpha}\hat{f}(\xi)|^{2}d\rho(\xi) \right\} ^{1/2}\\
  & = & \left\{ \sum _{|\alpha|=m}\frac{m!}{\alpha!}\int _{R^{n}}\xi ^{2\alpha}|\hat{f}(\xi)|^{2}d\rho(\xi)\right\} ^{1/2}\\
  & = & (m!)^{1/2}\left\{ \sum_{|\alpha|=m}\frac{1}{\alpha!}\int_{R^{n}}\xi^{2\alpha}|\hat{f}(\xi)|^{2}\frac{1}{(2\pi)^{2n}|\xi|^{2m}\hat{h}(\xi)}d\xi \right\} ^{1/2}(by\ (3.8)\ and\ Theorem\ 5.2\ of \cite{MN3})\\
  & = & \frac{\sqrt{m!}}{(2\pi)^{n}} \left\{ \sum_{|\alpha|=m}\frac{1}{\alpha!}\int_{R^{n}}\frac{\xi^{2\alpha}|\hat{f}(\xi)|^{2}}{|\xi|^{2m}}\cdot \frac{1}{2^{1+\frac{\beta}{2}}(\frac{|\xi|}{c})^{-\frac{\beta}{2}-\frac{n}{2}}{\cal K}_{\frac{n+\beta}{2}}(c|\xi|)}d\xi \right\} ^{1/2}(by\ Theorem\ 8.15\ of\ \cite{We})\\
  & \leq & \frac{\sqrt{m!}}{(2\pi)^{n}}\cdot 2^{-(\frac{1}{2}+\frac{\beta}{4})}\left\{ S(m,n)\int_{R^{n}}|\hat{f}(\xi)|^{2}\cdot \frac{|\xi|^{\frac{n+\beta}{2}}}{c^{\frac{n+\beta}{2}}}\cdot \frac{1}{{\cal K}_{\frac{n+\beta}{2}}(c|\xi|)}d\xi \right\} ^{1/2}\ where\ S(m,n)\ denotes\ the\\
  &   & number\ of\ terms\ in\ the\ ``\sum"\\
  & \leq & \frac{\sqrt{m!S(m,n)}}{(2\pi)^{n}2^{\frac{1}{2}+\frac{\beta}{4}}}\cdot c^{\frac{-(n+\beta)}{4}}\left\{ \int_{R^{n}}|\hat{f}(\xi)|^{2}|\xi|^{\frac{n+\beta}{2}}\cdot \frac{\sqrt{c|\xi|}}{\sqrt{\frac{\pi}{2}}e^{-c|\xi|}}d\xi \right\} ^{1/2}(by\ Corollary\ 5.12\ of\ \cite{We})\\
  & = & \frac{\sqrt{m!S(m,n)}}{(2\pi)^{n}2^{\frac{1}{2}+\frac{\beta}{4}}}\cdot c^{\frac{1-(n+\beta)}{4}}\left(\sqrt{\frac{2}{\pi}}\right)^{\frac{1}{2}}\left\{ \int_{R^{n}}|\hat{f}(\xi)|^{2}|\xi|^{\frac{1+n+\beta}{2}}e^{c|\xi|}d\xi \right\} ^{1/2}\\
  & \leq & \sqrt{m!S(m,n)}\cdot 2^{-n-\frac{1}{4}-\frac{\beta}{4}}\pi^{-\frac{1}{4}-n}\sigma^{\frac{1+n+\beta}{4}}e^{\frac{c\sigma}{2}}c^{\frac{1-(n+\beta)}{4}}\| f\| _{L^{2}(R^{n})} \\
  & < & \infty. 
\end{eqnarray*}
Thus $B_{\sigma}\subseteq {\cal C}_{h,m}$ and (7) follows.  \hspace{10cm} $\sharp$
\\

If in (1) $\beta<0$, the norm $\| \cdot \| _{h}$ is defined in a slightly different way. Hence we handle it separately and present the following theorem.
\begin{lem}
  Let $h$ be as in (1) with $\beta<0$ such that $n+\beta \geq 1$ or $n+\beta=-1$. Any function $f$ in $B_{\sigma}$ belongs to ${\cal C}_{h,m}$ and satisfies
\begin{equation}
  \| f\| _{h}\leq 2^{-n-\frac{1+\beta}{4}}\pi^{-n-\frac{1}{4}}\sigma^{\frac{1+\beta+n}{4}}e^{\frac{c\sigma}{2}}c^{\frac{1-(n+\beta)}{4}}\| f\| _{L^{2}(R^{n})}.
\end{equation}
\end{lem}
{\bf Proof}: If $\beta<0$, then $h$ is conditionally positive definite of order $m=0$, Wendland \cite{We}. By Theorem 5.2 of Madych et al. \cite{MN3}, we know that in Corollary 3.3 of \cite{MN3} $a_{\gamma}=0$ for $|\gamma|=2m$. Obviously any $f\in B_{\sigma}$ belongs to ${\cal S}'$ since $f\in L^{2}(R^{n})$. In order to apply Corollary 3.3 of \cite{MN3} to show that $B_{\sigma}\subseteq {\cal C}_{h,m}$, it remains to show that $f\in B_{\sigma}$ implies $\hat{f}\in L^{2}(\rho)$ where $d\rho (\xi):=r(\xi)d\xi =\frac{1}{(2\pi)^{2n}\hat{h}(\xi)}d\xi$ as stated in \cite{MN3}.

Now, let $f\in B_{\sigma}$. By (3.9) of \cite{MN3}, it suffices to show that $\| f\| _{h}:=\| \hat{f}\| _{L^{2}(\rho)}<\infty$. We proceed as follows.
\begin{eqnarray*}
  \| \hat{f}\| _{L^{2}(\rho)} & = & \left\{ \int_{R^{n}}|\hat{f}(\xi)|^{2}d\rho(\xi)\right\} ^{1/2}\\
                              & = & \left\{ \int_{R^{n}}|\hat{f}(\xi)|^{2}r(\xi)d\xi \right\} ^{1/2}\\
                              & = & \left\{ \int_{R^{n}}|\hat{f}(\xi)|^{2}(2\pi)^{-2n}\cdot \frac{1}{w(-\xi)}d\xi \right\} ^{1/2}\ where\ w(\cdot )=\hat{h}(\cdot )\ by\ Theorem\ 5.2\ of \ \cite{MN3}\\
                              & = & \left\{\int_{R^{n}}|\hat{f}(\xi)|^{2}(2\pi)^{-2n}2^{-1-\frac{\beta}{2}}(\frac{|\xi|}{c})^{\frac{\beta +n}{2}}\cdot \frac{1}{{\cal K}_{\frac{\beta +n}{2}}(c|\xi|)} d\xi \right\} ^{1/2}\ by \ (8.7),\ p.109\ of\  \cite{We} \\
                              & \leq & \left\{ 2^{-1-\frac{\beta}{2}-2n}\pi^{-2n}c^{-\frac{(n+\beta)}{2}}\int_{R^{n}}\frac{|\hat{f}(\xi)|^{2}|\xi|^{\frac{n+\beta}{2}}}{\sqrt{\frac{\pi}{2}}\cdot \frac{e^{-c|\xi|}}{\sqrt{c|\xi|}}}d\xi \right\} ^{1/2} \ if\ |n+\beta|\geq 1\ by\ Corollary\ 5.12\\
                              &   &  of\ \cite{We}\\
                              & = & \left\{ 2^{-\frac{1}{2}-\frac{\beta}{2}-2n}\pi^{-2n-\frac{1}{2}}c^{\frac{1}{2}-\frac{n+\beta}{2}}\int_{R^{n}}|\hat{f}(\xi)|^{2}|\xi|^{\frac{n+\beta+1}{2}}e^{c|\xi|}d\xi\right\} ^{1/2} \\
                              & \leq & \left\{ 2^{\frac{-1-\beta}{2}-2n}\pi^{-2n-\frac{1}{2}}c^{\frac{1}{2}-\frac{n+\beta}{2}}\sigma^{\frac{n+\beta+1}{2}}\int_{R^{n}}|\hat{f}(\xi)|^{2}e^{c|\xi|}d\xi\right\} ^{1/2}\ if\ n+\beta+1\geq 0 \\
                              & \leq & \left\{ 2^{\frac{-1-\beta}{2}-2n}\pi^{-2n-\frac{1}{2}}\sigma^{\frac{n+\beta+1}{2}}c^{\frac{1}{2}-\frac{n+\beta}{2}}e^{c\sigma}\int_{|\xi|\leq \sigma}|\hat{f}(\xi)|^{2}d\xi\right\} ^{1/2} \\
                              & \leq & 2^{\frac{-1-\beta}{4}-n}\pi^{-n-\frac{1}{4}}\sigma^{\frac{n+\beta+1}{4}}c^{\frac{1-(n+\beta)}{4}}e^{\frac{c\sigma}{2}}\| f\|_{L^{2}(R^{n})}\\                              & < & \infty \hspace{12.5cm} \sharp                             
\end{eqnarray*}

In the preceding deduction we put the restriction $|n+\beta|\geq 1$ on it. This is a drawback because the frequently seen case $n=1$ and $\beta=-1$ is not covered. This gap will be discussed shortly. Now Theorem 1.4 and Lemma 1.7 lead to the following theorem.
\begin{thm}
  For any positive number $\sigma$ and any function $f\in B_{\sigma}$, under the conditions of Theorem 1.4,
\begin{equation}
  |f(x)-s(x)|\leq \sqrt{m!S(m,n)}(2\pi)^{-\frac{3n}{4}}\sqrt{n\alpha_{n}}\sigma^{\frac{1+\beta+n}{4}}\sqrt{\Delta_{0}}c^{\frac{1+\beta-n}{4}}e^{\frac{c\sigma}{2}}(\lambda)^{\frac{1}{\delta}}\| f\|_{L^{2}(R^{n})}
\end{equation}
if $n+\beta\geq 1$ or $n+\beta=-1$. The number $S(m,n)$ is determined by $m$ and $n$, and is one whenever $\beta<0$. The other constants are as in Theorem 1.4.
\end{thm}

\begin{lem}
  Let $h$ be as in (1) with $\beta=-1,\ n=1$. For any $\sigma>0$, if $f\in B_{\sigma}$, then $f\in {\cal C}_{h,m}$ and 
\begin{equation}
\| f\| _{h}\leq (2\pi)^{-n}2^{-\frac{1}{4}}\left\{ \frac{1}{{\cal K}_{0}(1)}\int_{|\xi|\leq \frac{1}{c}}|\hat{f}(\xi)|^{2}d\xi + \frac{1}{a_{0}}\int_{\frac{1}{c}<|\xi|\leq \sigma}|\hat{f}(\xi)|^{2}\sqrt{c|\xi|}e^{c|\xi|}d\xi\right\} ^{1/2}
\end{equation}
if $\frac{1}{c}<\sigma$, where $a_{0}=\frac{\sqrt{\pi}3^{-\frac{1}{2}}}{2\Gamma(\frac{1}{2})}$, and
\begin{equation}
\| f\|_{h}\leq (2\pi)^{-n}2^{-\frac{1}{4}}\left\{ \frac{1}{{\cal K}_{0}(1)}\int_{|\xi|\leq \frac{1}{c}}|\hat{f}(\xi)|^{2}d\xi\right\} ^{1/2}
\end{equation}
if $\frac{1}{c}\geq \sigma $. Here ${\cal K}_{0}(\cdot )$ is the modified Bessel function. 
\end{lem}
Proof. By (3.9) of Madych et al. \cite{MN3}, $\|f\|_{h}=\|\hat{f}\|_{L^{2}(\rho)}$ and
\begin{eqnarray*}
  \|\hat{f}\|_{L^{2}(\rho)} & = & \left\{ \int_{R^{n}}|\hat{f}(\xi)|^{2}d\rho(\xi)\right\}^{1/2}\\
                            & = & \left\{ \int_{R^{n}}|\hat{f}(\xi)|^{2}r(\xi)d\xi\right\}^{1/2}\\
                            & = & \left\{ \int_{R^{n}}|\hat{f}(\xi)|^{2}(2\pi)^{-2n}\cdot \frac{1}{w(-\xi)}d\xi\right\}^{1/2}\ where\ w(\cdot)=\hat{h}(\cdot)\ by\ Theorem\ 5.2\ of\ \cite{MN3}\\
                            & = & \left\{ \int_{R^{n}}|\hat{f}(\xi)|^{2}(2\pi)^{-2n}2^{-\frac{1}{2}}\cdot \frac{1}{{\cal K}_{0}(c|\xi|)}d\xi \right\}^{1/2}\ by\ (8.7)\ of\ \cite{We}\ (p.109)\\
                            & = & [(2\pi)^{-2n}2^{-\frac{1}{2}}]^{\frac{1}{2}}\left\{ \int_{|\xi|\leq \frac{1}{c}}\frac{|\hat{f}(\xi)|^{2}}{{\cal K}_{0}(c|\xi|)}d\xi +\int_{|\xi|>\frac{1}{c}}\frac{|\hat{f}(\xi)|^{2}}{{\cal K}_{0}(c|\xi|)}d\xi\right\}^{1/2}\\
                            & \leq & (2\pi)^{-n}2^{-\frac{1}{4}}\left\{ \frac{1}{{\cal K}_{0}(1)}\int_{|\xi|\leq \frac{1}{c}}|\hat{f}(\xi)|^{2}d\xi+ \frac{1}{a_{0}}\int_{|\xi|>\frac{1}{c}}|\hat{f}(\xi)|^{2}\sqrt{c|\xi|}e^{c|\xi|}d\xi\right\}^{1/2}\\
                            &   & by\ Corollary\ 5.12\ of\ \cite{We}\ and\ p.\ 374\ of\ Abramowitz\ et\ al.\ \cite{Ab}\ where\ a_{0}=\frac{\sqrt{\pi}}{2\Gamma(\frac{1}{2})\sqrt{3}}. 
\end{eqnarray*}
Since $\hat{f}(\xi)=0$ for $|\xi|>\sigma$, our conclusion follows immediately. \hspace{5.3cm} $\sharp$
\begin{thm}
  For any positive number $\sigma$ and any function $f\in B_{\sigma}$, under the conditions of Theorem 1.4 with $\beta=-1$ and $n=1$,
\begin{equation}
  |f(x)-s(x)|\leq \sqrt{2\pi}\sqrt{\Delta_{0}}(\lambda)^{\frac{1}{\delta}}\left\{\frac{A+B}{c}\right\}^{1/2}
\end{equation}
where $A=\frac{1}{{\cal K}_{0}(1)}\int_{|\xi|\leq \frac{1}{c}}|\hat{f}(\xi)|^{2}d\xi$, $B=\frac{2\sqrt{3}\Gamma(\frac{1}{2})}{\sqrt{\pi}}\int_{\frac{1}{c}<|\xi|\leq \sigma}|\hat{f}(\xi)|^{2}\sqrt{c|\xi|}e^{c|\xi|}d\xi$ if $\frac{1}{c}\leq \sigma$ and $B=0$ if $\frac{1}{c}\geq \sigma$.
\end{thm}
Proof. This is just an immediate result of Theorem 1.4 and Lemma 1.9. \hspace{4cm} $\sharp$ 

\begin{cor}
  For any positive number $\sigma$ and any function $f\in B_{\sigma}$, under the conditions of Theorem 1.4 with $\beta=-1$ and $n=1$,
\begin{equation}
  |f(x)-s(x)|\leq \sqrt{2\pi}\sqrt{\Delta_{0}}(\lambda)^{\frac{1}{\delta}}M(c)\|f\|_{L^{2}(R^{n})}
\end{equation}
where $M(c):=\frac{1}{\sqrt{c}}\left \{\frac{1}{{\cal K}_{0}(1)}+\frac{2\sqrt{3}\Gamma(\frac{1}{2})}{\sqrt{\pi}}\sqrt{c\sigma}e^{c\sigma}\right \}^{\frac{1}{2}}$
\end{cor}

The proof is routine and we omit it.

\section{How to choose $c$?}

\subsection{$b_{0}$ fixed}
 We first deal with the situation when the domain $E$ is a cube of fixed side length $b_{0}$.

Note that in Theorem1.4
$$C=\max \left\{2\rho'\sqrt{n}e^{2n\gamma_{n}},\ \frac{2}{3b_{0}}\right\},\ \rho'=\frac{\rho}{c}.$$
Also, $2\rho'\sqrt{n}e^{2n\gamma_{n}}=\frac{2}{3b_{0}}$ if and only if $c=3b_{0}\rho\sqrt{n}e^{2n\gamma_{n}}$. Let $ c_{0}:=3b_{0}\rho\sqrt{n}e^{2n\gamma_{n}}$. In Theorem 1.4 it's required that $\delta\leq \delta_{0}$. For $c\leq c_{0}$, we have $\delta\leq \delta_{0}$ iff $c\geq c_{1}$ where $c_{1}:=12\rho \sqrt{n}e^{2n\gamma_{n}}\gamma_{n}(m+1)\delta$. Then $C=2\rho'\sqrt{n}e^{2n\gamma_{n}}$ if $c\in [c_{1},\ c_{0}]$, and $C=\frac{2}{3b_{0}}$ if $c\in [c_{0},\ \infty)$. For $c\in [c_{1},\ c_{0}]$,
\begin{eqnarray*}
  \lambda^{\frac{1}{\delta}} & = & \left(\frac{2}{3}\right)^{\frac{1}{6C\gamma_{n}\delta}}\\
                             & = & \left(\frac{2}{3}\right)^{\frac{c}{12\rho\sqrt{n}e^{2n\gamma_{n}}\gamma_{n}\delta}}\\
                             & = & \left[\left(\frac{2}{3}\right)^{\frac{1}{12\rho\sqrt{n}e^{2n\gamma_{n}}\gamma_{n}\delta}} \right]^{c}\\
                             & = & e^{\eta (\delta)c} 
\end{eqnarray*}
where $\eta(\delta)<0$ is some number depending on $\delta$. In fact $\eta(\delta):=\frac{ln\frac{2}{3}}{12\rho\sqrt{n}e^{2n\gamma_{n}}\gamma_{n}\delta}$. For $c\in [c_{0},\ \infty)$,
\begin{eqnarray*}
  \lambda^{\frac{1}{\delta}} & = & \left(\frac{2}{3}\right)^{\frac{1}{6C\gamma_{n}\delta}}\\
                             & = & \left(\frac{2}{3}\right)^{\frac{b_{0}}{4\gamma_{n}\delta}}.
\end{eqnarray*}

There are two cases.\\
\\

\hspace{-6mm}{\bf Case 1.} \fbox{$n+\beta\geq 1$ or $n+\beta=-1$}
 In (9), for $c\in [c_{1},\ c_{0}]$, the part influenced by $c$ is $c^{\frac{1+\beta-n}{4}}e^{\frac{c\sigma}{2}}(\lambda)^{\frac{1}{\delta}}=c^{\frac{1+\beta-n}{4}}e^{[\eta(\delta)+\frac{\sigma}{2}]c}$.

 Let's define $H_{1}(c):=c^{\frac{1+\beta-n}{4}}e^{[\eta(\delta)+\frac{\sigma}{2}]c}$ for $c\in [c_{1},\ c_{0}]$.

For $c\in [c_{0},\ \infty)$, we define $ H_{2}(c):=c^{\frac{1+\beta-n}{4}}e^{\frac{c\sigma}{2}}(\lambda)^{\frac{1}{\delta}}$ where $\lambda=(\frac{2}{3})^{\frac{1}{6C\gamma_{n}}}=(\frac{2}{3})^{\frac{b_{0}}{4\gamma_{n}}}$.

Note that $C\geq \frac{2}{3b_{0}}$ if $c\in [c_{1},\ c_{0}]$ and $C=\frac{2}{3b_{0}}$ if $c\in [c_{0},\ \infty)$.

Now let $$MN(c):= \left\{ \begin{array}{ll}
                          H_{1}(c) & \mbox{if $c\in [c_{1}, c_{0}]$}\\
                          H_{2}(c) & \mbox{if $c\in [c_{0},\ \infty)$.}
                        \end{array}\right. $$
Then minimizing the error bound (9) is equivalent to minimizing $MN(c)$.

The function $MN(c)$ is continuous and $H_{1}(c_{0})=H_{2}(c_{0})$. Theoretically, $c^{*}$ is the optimal choice of $c$ if $MN(c^{*})$ is the minimum of $MN(c)$.\\
\\
{\bf Examples:} In the following graphs $\delta$ denotes the fill distance of the interpolation centers defined in Theorem 1.4 and $b_{0}$ is the side length of the interpolation domain square in $R^{2}$.

\begin{figure}[h]
\centering
\includegraphics[scale=1.0]{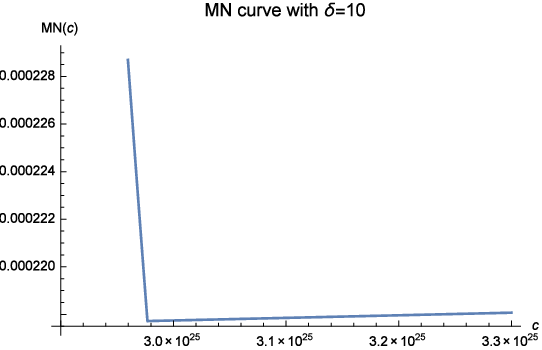}
\caption{Here $n=2,\ \beta=1, b_{0}=10000$ and $\sigma=10^{-27}$.}

\includegraphics[scale=1.0]{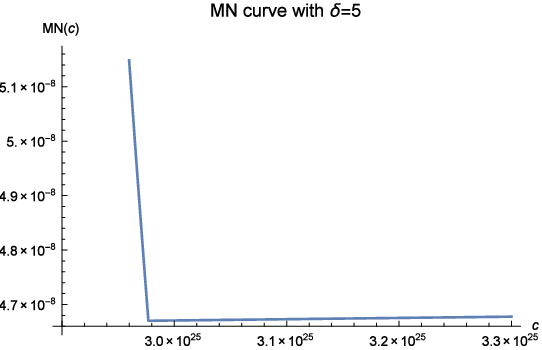}
\caption{Here $n=2,\ \beta=1, b_{0}=10000$ and $\sigma=10^{-27}$.}

\end{figure}

\clearpage

\begin{figure}[t]
\centering
\includegraphics[scale=1.0]{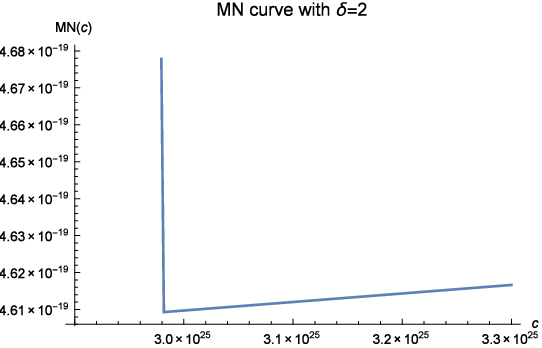}
\caption{Here $n=2,\ \beta=1, b_{0}=10000$ and $\sigma=10^{-27}$.}

\includegraphics[scale=1.0]{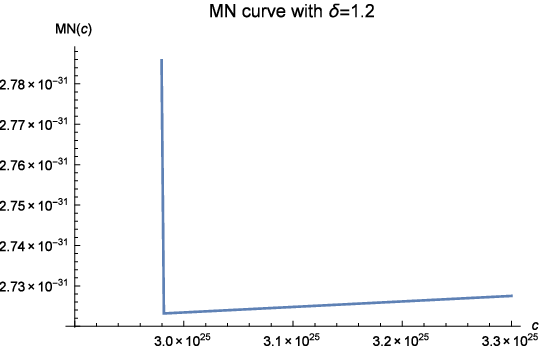}
\caption{Here $n=2,\ \beta=1, b_{0}=10000$ and $\sigma=10^{-27}$.}

\includegraphics[scale=1.0]{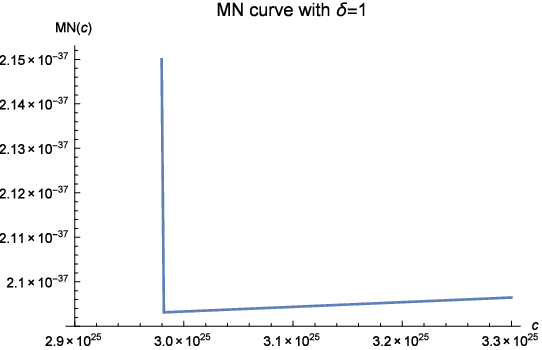}
\caption{Here $n=2,\ \beta=1, b_{0}=10000$ and $\sigma=10^{-27}$.}

\end{figure}

\clearpage

\hspace{-6mm}{\bf Case 2.} \fbox{$n=1$ and $\beta=-1$} In (13) the upper bound influenced by $c$ is$(\lambda)^{\frac{1}{\delta}}M(c)$ where $M(c)$ was defined in (13). For $c\in [c_{1},\ c_{0}]$, $\lambda^{\frac{1}{\delta}}=e^{\eta(\delta)c}$. For $c\in [c_{0},\ \infty)$, $\lambda^{\frac{1}{\delta}}=\left (\frac{2}{3}\right )^{\frac{b_{0}}{4\gamma_{n}\delta}}$. We thus define
 $$MN(c):= \left\{ \begin{array}{ll}
                          H_{1}(c) & \mbox{if $c\in [c_{1}, c_{0}]$,}\\
                          H_{2}(c) & \mbox{if $c\in [c_{0},\ \infty)$,}
                        \end{array}\right. $$
where $H_{1}:=e^{\eta(\delta)c}M(c)$ and $H_{2}:= \left(\frac{2}{3}\right)^{\frac{b_{0}}{4\gamma_{n}\delta}}M(c)$ with $M(c):=\frac{1}{\sqrt{c}}\left \{\frac{1}{{\cal K}_{0}(1)}+\frac{2\sqrt{3}\Gamma(\frac{1}{2})}{\sqrt{\pi}}\sqrt{c\sigma}e^{c\sigma}\right \}^{\frac{1}{2}}$ . The value $c$ minimizing $MN(c)$ is just the optimal value.\\
\\
{\bf Examples:} In the following graphs $\delta$ still denotes the fill distance, but $b_{0}$ is the length of the interpolation interval in $R^{1}$.
\begin{figure}[h]
\centering
\includegraphics[scale=1.0]{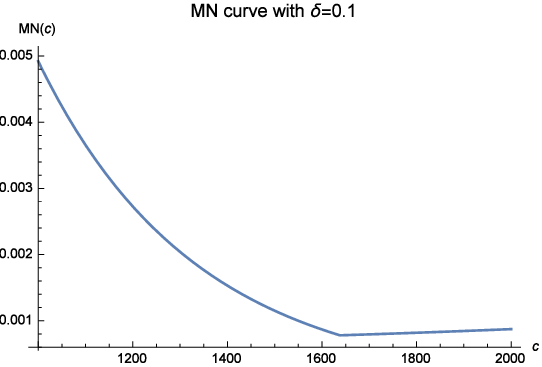}
\caption{Here $n=1,\ \beta=-1, b_{0}=10$ and $\sigma=10^{-3}$.}

\includegraphics[scale=1.0]{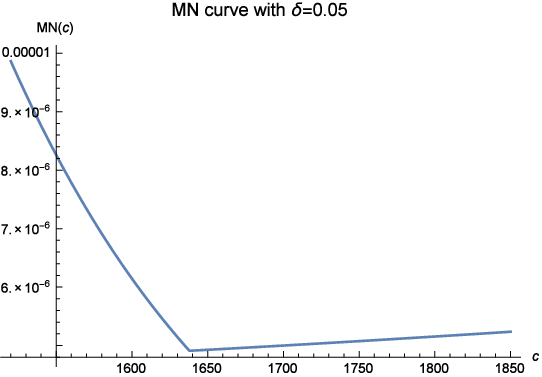}
\caption{Here $n=1,\ \beta=-1, b_{0}=10$ and $\sigma=10^{-3}$.}

\end{figure}

\clearpage

\begin{figure}[h]
\centering
\includegraphics[scale=0.9]{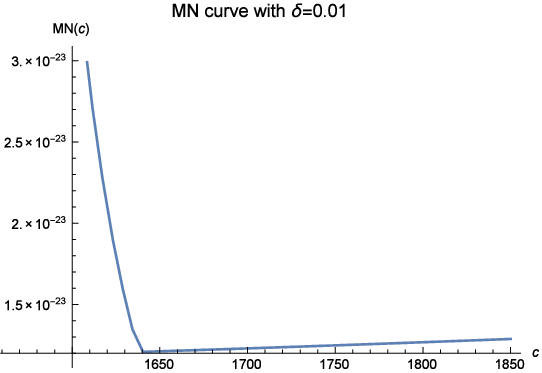}
\caption{Here $n=1,\ \beta=-1, b_{0}=10$ and $\sigma=10^{-3}$.}

\includegraphics[scale=0.9]{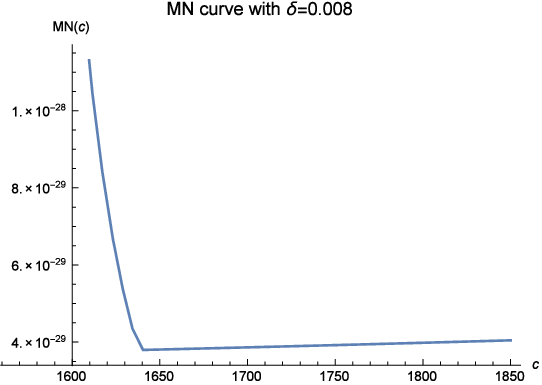}
\caption{Here $n=1,\ \beta=-1, b_{0}=10$ and $\sigma=10^{-3}$.}

\includegraphics[scale=0.9]{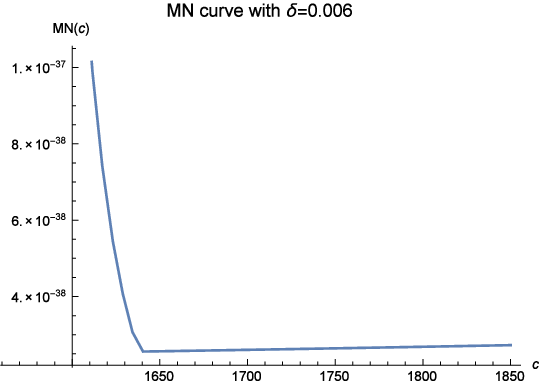}
\caption{Here $n=1,\ \beta=-1, b_{0}=10$ and $\sigma=10^{-3}$.}

\end{figure}

\clearpage

\hspace{-5.5mm}{\bf Remark:} One may wonder what will happen for $c\in (0,\ c_{1})$. Due to theoretical restriction, we didn't allow $c< c_{1}$. However, a huge number of experiments show that the optimal value of $c$ never falls into this interval. Hence we have essentially dealt with the entire interval $(0,\ \infty)$. Also, although the essential error bounds $MN(c)$ presented in Figs. 1-10 are already very small, the actual gap between the approximted function $f(x)$ and the approximating function $s(x)$ is much smaller, as we shall see in the next section.

\subsection{$b_{0}$ not fixed}
A lot of time the domain size can be determined by ourselves. If changing the side length $b_{0}$ of the interpolation domain $E$ makes the MN function value $MN(c)$ smaller, we should do it. In the error bounds (9) and (13), only $\lambda$ is influenced by $b_{0}$. In the structure of $\lambda$, if $C$ can be made smaller, $\lambda$ will become smaller, making the error bounds sharper. Note that decreasing $b_{0}$ makes $\lambda$ larger and is not worth doing. Therefore, without special reason, we will not decrease $b_{0}$. If both $c$ and $b_{0}$ increase, both $C$ and $\lambda$  will become smaller, making the error bounds better. There are two cases: (i) the fill distance $\delta$ is fixed, and (ii) $\delta$ is not fixed. For (i), if $b_{0}$ increases and $\delta$ is fixed, then the number of interpolation points will increases. This approach probably is only theoretically meaningful because it takes too much computer time. As for (ii),  if both $b_{0}$ and $\delta$ increase at the same speed, keeping $\frac{b_{0}}{\delta}$ constant in $\lambda^{\frac{1}{\delta}}$, then $MN(c)$ can be made smaller by increasing $c$, without increasing the number of the interpolation points. This is worth doing. In fact, some differential equations require that both the domain size and  distance among interpolation centers expand at the same speed as the time increases. For both (i) and (ii), we choose $c$ which minimizes the MN function, as in the preceding subsection.

\section{Experiment}
\subsection{\bf Small domain}
We first test Case 1 with $n=1,\ \beta=1,\ \sigma=10^{-4}$ and the interpolation domain $[0,\ b_{0}]$ where $b_{0}=10$.
In this experiment, the approximated function $f\in B_{\sigma}$ is $$f(x):=\frac{\sin(10^{-4}x)}{10^{-4}x}$$ where $f(x):=1$ if $x=0$. The approximating function is $s(x)$ as defined in (2). For simplicity, we replace the radial function $h(x)$ in (1) by the one presented in the Abstract. The interpolation centers are scattered in the interval $[0,\ 10]$. We use the Mathematica command RandomReal[$\cdot$] to generate a random real number $x_{i}'$ in the interval $[0,\ 1]$. Then let $x_{i}=\frac{10}{N_{d}}(i-1+ x_{i}')$, for $i=1, \ldots , N_{d}$. The interpolation centers $x_{i},\ldots , x_{N_{d}}$ will be scattered in $[0,\ 10]$ with fill distance $\delta =\frac{10}{N_{d}}$. In order to test the quality of the approximation, $N_{t}$ test points equally spaced in $[0,\ 10]$ are adopted. The quality of the approximation is evaluated by the root-mean-square error defined by
$$RMS:=\sqrt{\frac{\sum_{i=1}^{N_{t}}|f(\overline{x}_{i})-s(\overline{x}_{i})|^{2}}{N_{t}}}$$ 
where $\overline{x}_{1}, \ldots ,\overline{x}_{N_{t}}$ are the test points. The condition number of the interpolation linear system constituted by (3), (4) is denoted by $COND$.

Before entering our experiment, let us analyze a few MN curves.

\begin{figure}[t]
\centering
\includegraphics[scale=0.9]{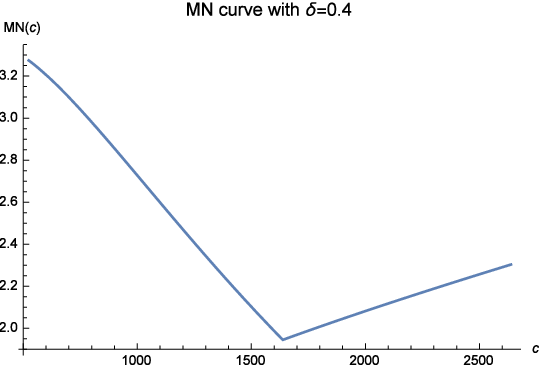}
\caption{Here $n=1,\ \beta=1, b_{0}=10$ and $\sigma=10^{-4}$.}

\includegraphics[scale=0.9]{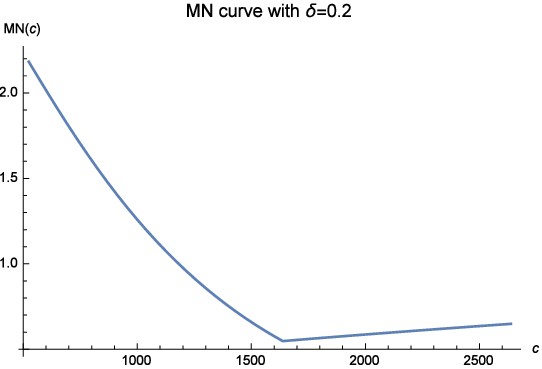}
\caption{Here $n=1,\ \beta=1, b_{0}=10$ and $\sigma=10^{-4}$.}

\includegraphics[scale=0.9]{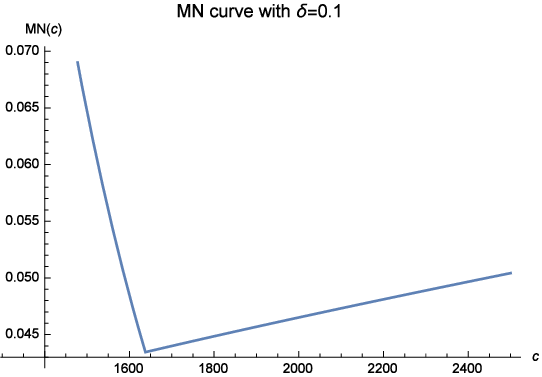}
\caption{Here $n=1,\ \beta=1, b_{0}=10$ and $\sigma=10^{-4}$.}

\end{figure}

\clearpage

\begin{figure}[t]
\centering
\includegraphics[scale=1.0]{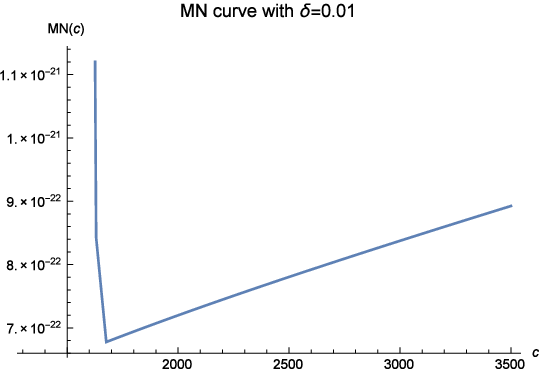}
\caption{Here $n=1,\ \beta=1, b_{0}=10$ and $\sigma=10^{-4}$.}
\end{figure}

Figures 11-14 strongly suggest that one should choose $c$ to be the number $c_{0}=1637.94$. We test $\delta=0.4,\ 0.2$, and 0.1. The results are presented in Tables 1-3, where $c_{0}=1637.94$ and $c_{1}=1408.28,\ 524.124,\ 262.071$, respectively. In order to make the tables in this paper easier to understand, we offer a summary  for the definitions of the notations as follows.\\
\\
{\bf Notations:} $\delta$: the fill distance; $b_{0}$: the side length of the domain cube; $N_{d}$: the number of interpolation points used; $N_{t}$: the number of test points; $c$: the shape parameter contained in the radial function; $c_{1}$: the smallest $c$ allowed by our theory; $c_{0}$: the joint point of $H_{1}(c)$ and $H_{2}(c)$ in $MN(c)$; $COND$: the condition number of the interpolation linear system; $RMS$: the root-mean-square error of the interpolation

\begin{table}[h]
\caption{$\delta=0.4,\ b_{0}=10,\ N_{d}=25,\ N_{t}=1000$}
\centering
\tiny
\begin{tabular}{c lllll}\\[2ex]
\hline\hline \\ [1ex]
\large $c$ & \large 1 & \large 1000  & \large $c_{1}$  & \large $c_{0}$ & \large 5000  \\ [1ex]
\hline\\[1ex]

\large $RMS$ & \normalsize $ 5.16\cdot 10^{-11} $ & \normalsize $3.74\cdot 10^{-64}$ & \normalsize $8.74\cdot 10^{-65} $   & \normalsize $ 2.67\cdot 10^{-69} $ & \normalsize $5.74\cdot 10^{-81} $     \\ [1ex]
\hline\\[1ex]

\large $COND$ & \normalsize $2.08\cdot 10^{6}$ & \normalsize $2.07\cdot 10^{124}$ & \normalsize $7.13\cdot 10^{124}$   & \normalsize $4.01 \cdot 10^{134}$ & \normalsize $7.36\cdot 10^{157}$  \\ [1ex]

\hline\hline \\[1ex]

\large $c$ & \large $10^{4}$ & \large $10^{5}$   \\ [1ex]
\hline\\[1ex]

\large $RMS$ & \normalsize $ 2.52\cdot 10^{-88} $ & \normalsize $4.94\cdot 10^{-112} $     \\ [1ex]
\hline \\[1ex]

\large $COND$ & \normalsize $2.07\cdot 10^{172}$ & \normalsize $1.25\cdot 10^{220}$  \\ [1ex]
\hline \\[1ex]

\end{tabular}

\label{001}
\end{table}

\clearpage

\begin{table}[h]
\caption{$\delta=0.2,\ b_{0}=10,\ N_{d}=50,\ N_{t}=1000$}
\centering
\tiny
\begin{tabular}{c lllll}\\[2ex]
\hline\hline \\ [1ex]
\large $c$ & \large 1 & \large $c_{1}$ & \large 1000 & \large $c_{0}$ & \large 5000  \\ [1ex]
\hline\\[1ex]

\large $RMS$ & \normalsize $ 1.14\cdot 10^{-12} $ & \normalsize $4.76\cdot 10^{-105} $ & \normalsize $ 1.63\cdot 10^{-118} $ & \normalsize $ 8.38\cdot 10^{-129} $ & \normalsize $4.28\cdot 10^{-152} $     \\ [1ex]
\hline\\[1ex]

\large $COND$ & \normalsize $3.62\cdot 10^{10}$ & \normalsize $2.70\cdot 10^{224}$ & \normalsize $8.42 \cdot 10^{251}$ & \normalsize $8.44 \cdot 10^{272}$ & \normalsize $2.65\cdot 10^{320}$  \\ [1ex]

\hline\hline \\[1ex]

\large $c$ & \large $10^{4}$ & \large $5\cdot 10^{6}$   \\ [1ex]
\hline\\[1ex]

\large $RMS$ & \normalsize $3.07\cdot 10^{-166}$ & \normalsize $4.65\cdot 10^{-240} $     \\ [1ex]
\hline \\[1ex]

\large $COND$ & \normalsize $7.36\cdot 10^{349}$ & \normalsize $2.65\cdot 10^{614}$  \\ [1ex]
\hline \\[1ex]

\end{tabular}

\label{002}
\end{table}

\begin{table}[h]
\caption{$\delta=0.1,\ b_{0}=10,\ N_{d}=100,\ N_{t}=1000$}
\centering
\tiny
\begin{tabular}{c lllll}\\[2ex]
\hline\hline \\ [1ex]
\large $c$ & \large 1 & \large 100 & \large$c_{1}$ & \large $1000$ & \large$c_{0}$  \\ [1ex]
\hline\\[1ex]

\large $RMS$ & \normalsize $ 5.05\cdot 10^{-14} $ & \normalsize $1.49\cdot 10^{-132} $ & \normalsize $ 9.13\cdot 10^{-174} $ & \normalsize $ 1.93\cdot 10^{-220} $ & \normalsize $1.92\cdot 10^{-251} $     \\ [1ex]
\hline\\[1ex]

\large $COND$ & \normalsize $3.91\cdot 10^{17}$ & \normalsize $1.00\cdot 10^{308}$ & \normalsize $4.03 \cdot 10^{390}$ & \normalsize $8.71 \cdot 10^{505}$ & \normalsize $2.35\cdot 10^{548}$  \\ [1ex]

\hline\hline \\[1ex]

\large $c$ & \large $5000$ & \large $10^{4}$   \\ [1ex]
\hline\\[1ex]

\large $RMS$ & \normalsize $ 5.61\cdot 10^{-299} $ & \normalsize $1.30\cdot 10^{-328} $     \\ [1ex]
\hline \\[1ex]

\large $COND$ & \normalsize $2.17\cdot 10^{644}$ & \normalsize $8.70\cdot 10^{703}$  \\ [1ex]
\hline \\[1ex]

\end{tabular}

\label{003}
\end{table}

\subsection{\bf Large domain}
We now test Case 1 for $n=1,\ \beta=1,\ \sigma=10^{-33}$, and the interpolation domain $[0,\ b_{0}]$ where $b_{0}=10^{30}$. The approximated function $f\in B_{\sigma}$ is 
$$f(x):=\frac{\sin (10^{-33}x)}{10^{-33}x}$$ where $f(x):=1$ if $x=0$.
A few MN curves are offered as follows.

\begin{figure}[t]
\centering
\includegraphics[scale=0.9]{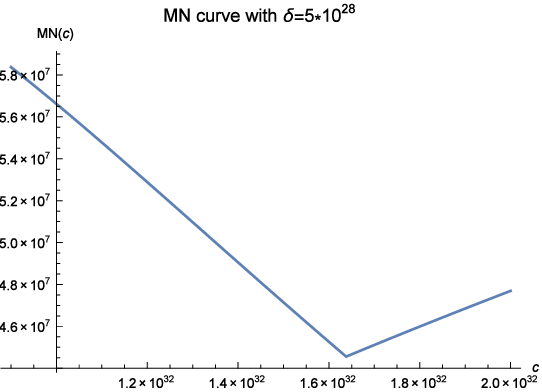}
\caption{Here $n=1,\ \beta=1, b_{0}=10^{30}$ and $\sigma=10^{-33}$.}

\includegraphics[scale=0.9]{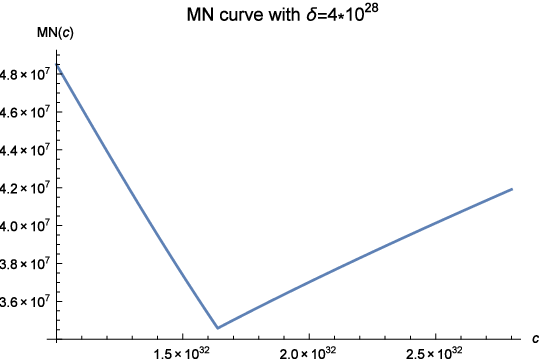}
\caption{Here $n=1,\ \beta=1, b_{0}=10^{30}$ and $\sigma=10^{-33}$.}

\includegraphics[scale=0.9]{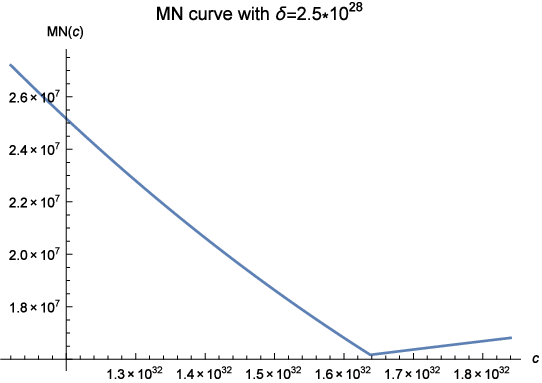}
\caption{Here $n=1,\ \beta=1, b_{0}=10^{30}$ and $\sigma=10^{-33}$.}

\end{figure}

\clearpage

\begin{figure}[t]
\centering
\includegraphics[scale=1.1]{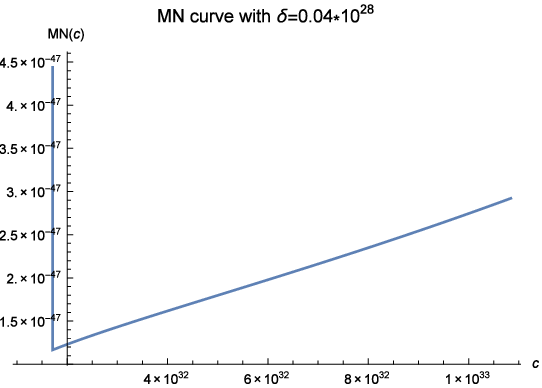}
\caption{Here $n=1,\ \beta=1, b_{0}=10^{30}$ and $\sigma=10^{-33}$.}

\end{figure}

We test $\delta=5\cdot 10^{28}$, $\delta=4\cdot 10^{28}$ and $\delta=2.5\cdot 10^{28}$. The results are shown in Tables 4-6, where $c_{0}=1.63794\cdot 10^{32}$ is the optimal choice of $c$ suggested by the MN curves, and $c_{1}$ is the smallest value of $c$ permitted by our theory. In Tables 4, 5, and 6, $c_{1}=1.31036\cdot 10^{32},\ 1.04828\cdot 10^{32}$ and $6.55178\cdot 10^{31}$, respectively.

\begin{table}[h]
\caption{$\delta=5\cdot 10^{28},\ b_{0}=10^{30},\ N_{d}=20,\ N_{t}=200$}
\centering
\tiny
\begin{tabular}{c lllll}\\[2ex]
\hline\hline \\ [1ex]
\large $c$ & \large $10^{24}$ & \large $10^{26}$ & \large $10^{28}$ & \large $10^{30}$ & \large $c_{1}$  \\ [1ex]
\hline\\[1ex]

\large $RMS$ & \normalsize $ 4.55\cdot 10^{-10} $ & \normalsize $4.5\cdot 10^{-10} $ & \normalsize $ 3.81\cdot 10^{-10} $ & \normalsize $ 3.61\cdot 10^{-16} $ & \normalsize $1.37\cdot 10^{-52} $     \\ [1ex]
\hline\\[1ex]

\large $COND$ & \normalsize $4.56\cdot 10^{60}$ & \normalsize $4.6\cdot 10^{60}$ & \normalsize $4.61 \cdot 10^{60}$ & \normalsize $4.62 \cdot 10^{61}$ & \normalsize $1.12\cdot 10^{103}$  \\ [1ex]

\hline\hline \\[1ex]

\large $c$ & \large $c_{0}$ & \large $10^{40}$ & \large $10^{50}$  \\ [1ex]
\hline\\[1ex]

\large $RMS$ & \normalsize $ 2.47\cdot 10^{-54} $ & \normalsize $6.03\cdot 10^{-90} $  & \normalsize $6.03\cdot 10^{-90}$   \\ [1ex]
\hline \\[1ex]

\large $COND$ & \normalsize $5.39\cdot 10^{106}$ & \normalsize $1.32\cdot 10^{402}$ & \normalsize $1.32\cdot 10^{782}$ \\ [1ex]
\hline \\[1ex]

\end{tabular}

\label{004}
\end{table}

\clearpage

\begin{table}[t]
\caption{$\delta=4\cdot 10^{28},\ b_{0}=10^{30},\ N_{d}=25,\ N_{t}=200$}
\centering
\tiny
\begin{tabular}{c lllll}\\[2ex]
\hline\hline \\ [1ex]
\large $c$ & \large $10^{20}$ & \large $10^{26}$ & \large $10^{28}$ & \large $10^{30}$ & \large $c_{1}$  \\ [1ex]
\hline\\[1ex]

\large $RMS$ & \normalsize $ 1.52\cdot 10^{-9} $ & \normalsize $1.52\cdot 10^{-9} $ & \normalsize $ 1.26\cdot 10^{-9} $ & \normalsize $ 1.05\cdot 10^{-17} $ & \normalsize $4.49\cdot 10^{-64} $     \\ [1ex]
\hline\\[1ex]

\large $COND$ & \normalsize $5.61\cdot 10^{60}$ & \normalsize $5.61\cdot 10^{60}$ & \normalsize $5.68 \cdot 10^{60}$ & \normalsize $6.0 \cdot 10^{61}$ & \normalsize $2.62\cdot 10^{125}$  \\ [1ex]

\hline\hline \\[1ex]

\large $c$ & \large $c_{0}$ & \large $10^{40}$   \\ [1ex]
\hline\\[1ex]

\large $RMS$ & \normalsize $ 9.95\cdot 10^{-69} $ & \normalsize $4.21\cdot 10^{-117} $    \\ [1ex]
\hline \\[1ex]

\large $COND$ & \normalsize $5.27\cdot 10^{134}$ & \normalsize $4.64\cdot 10^{508}$ \\ [1ex]
\hline \\[1ex]

\end{tabular}

\label{005}
\end{table}

\begin{table}[h]
\caption{$\delta=2.5\cdot 10^{28},\ b_{0}=10^{30},\ N_{d}=40,\ N_{t}=200$}
\centering
\tiny
\begin{tabular}{c lllll}\\[2ex]
\hline\hline \\ [1ex]
\large $c$ & \large $10^{26}$ & \large $10^{28}$ & \large $c_{1}$ & \large $c_{0}$& \large $10^{40}$  \\ [1ex]
\hline\\[1ex]

\large $RMS$ & \normalsize $ 6.01\cdot 10^{-10} $ & \normalsize $4.43\cdot 10^{-10} $ & \normalsize $ 1.38\cdot 10^{-88} $ & \normalsize $ 1.04\cdot 10^{-103} $ & \normalsize $1.18\cdot 10^{-188} $     \\ [1ex]
\hline\\[1ex]

\large $COND$ & \normalsize $9.75\cdot 10^{60}$ & \normalsize $9.87\cdot 10^{60}$ & \normalsize $2.56 \cdot 10^{186}$ & \normalsize $2.80 \cdot 10^{217}$ & \normalsize $6.92\cdot 10^{824}$  \\ [1ex]

\hline\hline \\[1ex]

\large $c$ & \large $ 10^{50}$  \\ [1ex]
\hline\\[1ex]

\large $RMS$ & \normalsize $ 3.42\cdot 10^{-188} $    \\ [1ex]
\hline \\[1ex]

\large $COND$ & \normalsize $6.92\cdot 10^{824}$  \\ [1ex]
\hline \\[1ex]

\end{tabular}

\label{006}
\end{table}

It is clearly seen that in Tables 1-6, the interpolation error $RMS$ is always very small when $c$ is equal to the optimal value $c_{0}$ suggested by the MN curves, even if this choice may not be the experimentally best one. In our previous papers Luh \cite{Lu4, Lu5, Lu6}, one can always find the experimentally best $c$ according to the MN curves by decreasing the value of $\delta$. However, here we scrapped this approach for two reasons. First, the $RMS$ is already very small if $c$ is chosen according to the MN curves, even when $\delta$ is quite large. Second, if we decrease the fill distance $\delta$, the condition numbers will become very large. In order to overcome ill-conditioning, one has to increase the number of effective digits for each step of the calculation, resulting in a huge amount of computer time. In fact, the results presented in Tables 1-6 are only a very small subset of our results. We of course have tested much smaller $\delta$'s and obtained better $c$'s, but the considerable computer time made us decide not to go further because it might not be worth doing. Note that, comparing Tables 1-6 and Figures 11-18, one will find that there is a big gap between the MN function value and the interpolation error $RMS$. It shows that the essential error bound $MN(c)$ is far from being sharp. Therefore, it does not seem to be sensible to expect for the experimentally optimal choice of $c$ merely by the MN curves. This is the main difference between our evenly-spaced-data and scattered-data schemes.The former is based on a much sharper error bound, and hence the theoretically and experimentally optimal values of $c$ always coincide. Fortunately, our scattered-data approach also provides very good results as long as $c$ is chosen according to the MN curves.

A practical question to be explained is how we overcame ill-conditioning in the experiment. By the help of the infinitely precise computer software Mathematica, we always adopted enough effective digits for each step of the calculation. For example, if $COND=10^{824}$, then we kept 1200 effective digits to the right of the decimal point. If increasing the number of effective digits resulted in the same $RMS$, we stopped it. This approach works well. However, if $c<<\|x_{i}-x_{j}\|$, then $\sqrt{c^{2}+\|x_{i}-x_{j}\|^{2}}\approx \|x_{i}-x_{j}\|$, making the interpolation matrix numerically remain unchanged by changing $c$ and hence the approximating function $s(x)$ also unchanged. This phenomenon can be seen in Table 5 for $c=10^{20}$ and $10^{26}$, even if enough effective digits were adopted. Conversely, if $c>>\|x_{i}-x_{j}\|$, then $\sqrt{c^{2}+\|x_{i}-x_{j}\|^{2}}\approx c$, and the final approximating function $s(x)$ obtained is numerically unchanged, even if the linear system is solvable by keeping enough effective digits. This can be seen in Table 4 for $c=10^{40}$ and $10^{50}$.

What is exciting is that very small $RMS$ can always be obtained by using very few data points, even if the domain is a huge interval. This can be seen in Tables 4-6, where the interpolation domain is $[0,\ 10^{30}]$, but only 20, 25, and 40 data points are involved.

\section{Final remarks}
The choice of the shape parameter has obsessed experts in the field of RBF since 1971. This paper unveils its mystery theoretically and presents a practically useful approach to choosing it. In the past decades, various error bounds for RBF interpolations were established, but only the one raised by Madych and Nelson in Madych et al. \cite{MN3} , which is just the cornerstone of our approach, leads to an essentially meaningful and useful criterion of choosing $c$ for purely scattered data settings. It is tempting to attempt to improve our aproach and make the essential error bound $MN(c)$ sharper. However, if one investigates the proofs in \cite{MN3} carefully, it can be seen that there is little opportunity for the improvement.

Finally, we will emphasize again on the applications of our approach to PDEs. Although the function space $B_{\sigma}$ in this paper is quite small, it plays just an intermediate role between the approximating function $s(x)$ defined in (2) and the approximated functions in the Sobolev space where the solutions of many important differential equations exist. Any function in the Sobolev space can be interpolated by a $B_{\sigma}$ function, with a good error bound, as shown in Narcowich et al. \cite{NWW}. Then any $B_{\sigma}$ function can be interpolated by $s(x)$, which is constructed by RBF, also with a good error bound. The final error estimates can be handled by triangle inequality. The $B_{\sigma}$ function need not be found explicitly. One only needs its existence and the two steps of interpolation share the same set of data points.

% Set the ending of a LaTeX document

\end{document}